\documentclass[12pt]{amsart}
\usepackage[usenames,dvipsnames]{xcolor}
\usepackage{latexsym}
\usepackage{amssymb}
\usepackage{amsmath}
\usepackage{amsthm}
\usepackage[T1]{fontenc}
\usepackage[latin2,cp1250]{inputenc}
\usepackage{indentfirst}
\usepackage{subcaption}
\usepackage{relsize}

\newtheorem{tw}{Theorem}
\newtheorem{pr}[tw]{Proposition}
\newtheorem{lem}[tw]{Lemma}
\newtheorem{co}[tw]{Corollary}
\newtheorem{fkt}[tw]{Fact}

\theoremstyle{definition}
\newtheorem{df}[tw]{Definition}
\newtheorem{re}[tw]{Remark}

\newcommand{\R}{\mathbb{R}}
\newcommand{\N}{\mathbb{N}}
\newcommand{\mB}{\mathcal{B}}
\newcommand{\mP}{\mathcal{P}}
\newcommand{\mN}{\mathcal{N}}
\newcommand{\mA}{\mathcal{A}}
\newcommand{\mM}{\mathcal{M}}
\newcommand{\xD}{(X,d)}
\newcommand{\xBF}{(X,\mB^\varphi)}
\newcommand{\st}{\mid}
\newcommand{\pruf}{\begin{proof}}
\newcommand{\pruuf}{\end{proof}}
\newcommand{\forol}[1]{\displaystyle\mathlarger{\mathop{\mathlarger{\mathlarger{\forall}}}_{#1}}\ }
\newcommand{\egzisc}[1]{\displaystyle\mathlarger{\mathop{\mathlarger{\mathlarger{\exists}}}_{#1}}\ }
\newcommand{\tvarphi}{{\tilde{\varphi}}}

\usepackage{hyperref}
\hypersetup{
	pdfstartview={FitB},
	colorlinks=true,
	linkcolor=Fuchsia,
	citecolor=SeaGreen
}

\begin{document}
\title[Caristi--Kirk and Oettli--Th\'era Ball Spaces]{Caristi--Kirk and Oettli--Th\'era Ball Spaces and applications}

\author[B\l aszkiewicz, \'Cmiel, Linzi, Szewczyk]{Piotr B\l aszkiewicz, Hanna \'Cmiel,\\
Alessandro Linzi, Piotr Szewczyk}
\address{Institute of Mathematics, University of Szczecin\newline ul.\ Wielkopolska 15, 70-451 Szczecin, Poland}
\email{blaszko333@gmail.com, hannacmielmath@gmail.com,\newline linzi.alessandro@gmail.com, piord.szewczyk@gmail.com}
\date{}
\begin{abstract}
Based on the theory of ball spaces introduced by Kuhlmann and Kuhlmann we introduce and study Caristi--Kirk and Oettli--Th\'era ball spaces. We show that if the underlying metric space is complete, then these have a very strong property: every ball contains a singleton ball. This fact provides quick proofs for several results which are equivalent to the Caristi--Kirk Fixed Point Theorem, namely Ekeland's Variational Principles, the Oettli--Th\'era Theorem, Takahashi's Theorem and the Flower Petal Theorem.

\end{abstract}
\maketitle
\section{Introduction}
\subsection{General setting}

The literature on complete metric spaces contains remarkable results such as the Theorem of Caristi and Kirk (\cite{caristi} and \cite{kirk}), Ekeland's Principle (\cite{ekeland}), Takahashi's Theorem (\cite{takahashi}) and the Flower Petal Theorem (\cite{penot}).
These theorems are known to be equivalent (see, e.g., \cite{penot}, \cite{ot}). Their statements can be found in Section \ref{original}.\par

The concept of a ball space was first introduced by F.-V. and K. Kuhlmann in \cite{kbonus},\cite{kbonus2}. In \cite{karticle} they connected it with the Caristi--Kirk Fixed Point Theorem (FPT) by providing a way to prove it using ball spaces techniques. In this paper we further develop this connection by proving Theorem \ref{poetycko} which provides a generic method to obtain simple proofs of all the results mentioned in the previous paragraph and, possibly, related ones in the future.\par
In \cite{ot}, Oettli and Th\'era introduced an alternative approach to the Caristi--Kirk Theorem and showed it to be equivalent to what was later (in publications such as \cite{meghea}) called Oettli--Th\'era Theorem. Our method can be applied to easily prove this theorem as well as the theorems equivalent to it, which are stated in \cite{ot} (see Section \ref{applications}). 

\subsection{Ball spaces}
As in \cite{karticle}, by a \textit{ball space} we mean a pair $(X,\mathcal{B})$, where $X$ is a nonempty set and $\mathcal{B}\subseteq\mP(X)$  is a nonempty family of nonempty subsets of $X$. An element $B\in \mathcal{B}$ is called a \textit{ball}. If no confusion arises, we will write $\mB$ in place of $(X,\mB)$ when speaking of a ball space.\par
A \textit{nest of balls} in a ball space $\mB$ is a nonempty family $\mathcal{N}$ of balls from $\mathcal{B}$ which is totally ordered by inclusion. We say that a ball space $\mB$ is \textit{spherically complete} if for every nest of balls $\mathcal{N}\subseteq\mathcal{B}$ we have $\bigcap\mathcal{N}\neq\emptyset$. Further details about ball spaces may be found in \cite{kbook}.
\begin{df}\label{stronc}
A ball space $(X,\mB)$ is \emph{strongly contractive} if there is a function that associates to every $x\in X$ some ball $B_x\in\mB$ such that, for every $x,y\in X$, the following conditions hold:
\begin{enumerate}
\item $x\in B_x$;
\item if $y\in B_x$ then $B_y\subseteq B_x$; 
\item if $y\in B_x\setminus\{x\}$ then $B_y\subsetneq B_x$.
\end{enumerate}
\end{df}
This particular type of ball spaces has a remarkable property, stated in the following theorem.
\begin{tw}
\label{poetycko}
In every spherically complete, strongly contractive ball space every ball $B_x$ contains a singleton ball. In other words, there exists $a\in B_x$ such that $B_a=\{a\}$.
\end{tw}
\pruf
Let $\mB$ be a strongly contractive, spherically complete ball space and $B_x\in\mB$ any ball. Consider the family
\[
\mA=\{\mN\subseteq\mP(B_x)\st\mN\text{ is a nest of balls in }\mB\}.
\]
This family is partially ordered by inclusion and nonempty since $\{B_x\}\in\mA$. If we have a chain of nests in $\mA$, the union of that chain is again a nest of balls in $\mA$, hence an upper bound of the chain. By Zorn's Lemma we obtain the existence of a maximal nest $\mM\in\mA$. Since the space is spherically complete, there exists an element $a\in\bigcap\mM$. Since $a\in B$ for every $B\in\mM$, by condition (2) of Definition \ref{stronc} also $B_a\subseteq B$ for every $B\in\mM$ and so $B_a\subseteq\bigcap\mM$. This means that $\mM\cup\{B_a\}$ is a nest of balls in $\mA$ which contains $\mM$. By maximality of $\mM$ we get that $\mM\cup\{B_a\}=\mM$, i.e., $B_a\in\mM$. Now we wish to show that $B_a$ is a singleton. Suppose that there exists an element $b\in B_a\setminus\{a\}$. Then $B_b\subsetneq B_a$ (in particular, $B_a\not\subseteq B_b$) and so $B_b\notin\mM$. But this means that $\mM\cup\{B_b\}$ is a nest of balls that properly contains $\mM$, which contradicts the maximality of $\mM$. Therefore, $B_a=\{a\}$.
\pruuf
\section{Caristi--Kirk and Oettli--Th\'era ball spaces}
In this section, we will be working with a nonempty metric space $\xD$.
\subsection{Caristi--Kirk ball spaces}
\label{CKBS}
Consider a function $\varphi:X\to\R$, a point $x\in X$ and the following set:
\[
B_x^\varphi=\{y\in X\st d(x,y)\le\varphi(x)-\varphi(y)\}.
\]
Since $B_x^\varphi\ne\emptyset$ (because $x\in B^\varphi_x$), we may think of this set as a ball and consider the ball space $(X,\mB^\varphi)$ where
\[
\mB^\varphi:=\left\{B_x^\varphi\st x\in X\right\}.
\]
We will call the function $\varphi$ a \emph{Caristi--Kirk function on }$X$ if it is lower semicontinuous, that is,
\[
\forol{y\in X}\liminf_{x\to y}\varphi(x)\ge\varphi(y),
\]
and bounded from below, that is,
\[
\inf_{x\in X}\varphi(x)>-\infty.
\]
The corresponding balls $B^\varphi_x$ have been introduced in \cite{karticle} as the \emph{Caristi--Kirk balls} and $\mB^\varphi$ is the induced \emph{Caristi--Kirk ball space}. For brevity, we will write CK in place of Caristi--Kirk.\par
A number of remarkable properties of the balls defined above, given in the following lemma, can be found in \cite{karticle}.
\begin{lem}
\label{lem0}
Take a metric space $\xD$ and any function $\varphi:X\to\R$. Then the following assertions hold.
\begin{enumerate}
\item For every $x\in X$, $x\in B^\varphi_x$.
\item If $y\in B^\varphi_x$ then $B^\varphi_y\subseteq B^\varphi_x$; if in addition $x\ne y$, then $B^\varphi_y\subsetneq B^\varphi_x$ and $\varphi(y)<\varphi(x)$.
\item If $\varphi$ is lower semicontinuous, then all CK balls $B_x$ are closed in the topology induced by the metric.
\end{enumerate}
\end{lem}
Lemma \ref{lem0} immediately yields the following result.
\begin{co}\label{co00}
The CK ball space $\mB^\varphi$ is strongly contractive.
\end{co}
Another important fact about CK ball spaces may also be found in \cite{karticle}:
\begin{pr}
\label{pr1}
Let $\xD$ be a metric space. Then the following statements are equivalent:
\begin{itemize}
\item[$(i)$] The metric space $\xD$ is complete.
\item[$(ii)$] Every CK ball space $\xBF$ is spherically complete.
\item[$(iii)$] For every continuous function $\varphi:X\to\R$ bounded from below, the CK ball space $\xBF$ is spherically complete.
\end{itemize}
\end{pr}

\subsection{Oettli--Th\'era ball spaces}
\label{OTBS}
\begin{df}\label{whoa}
A function $\phi:X\times X\to(-\infty,+\infty]$ is an \emph{Oettli--Th\'era function on }$X$ if the following properties hold:
\begin{align*}
(a)\quad&\phi(x,\cdot):X\rightarrow(-\infty,+\infty]\text{ is lower semicontinous for all } x\in X;\\
(b)\quad&\phi(x,x)=0\text{ for all }x\in X;\\
(c)\quad&\phi(x,y)\le \phi(x,z)+\phi(z,y)\text{ for all }x,y,z\in X;\\
(d)\quad&\text{there exists }x_0\in X\text{ s.t.\ } \inf_{x\in X}\phi(x_0,x)>-\infty.
\end{align*}
If an element $x_0\in X$ satisfies property $(d)$, we will call it an \emph{Oettli--Th\'era element for }$\phi$\emph{ in }$X$. If it is clear which space is considered, we will say that $\phi$ is an Oettli--Th\'era function and that $x_0$ is an Oettli--Th\'era element for $\phi$. For brevity, we will write OT in place of Oettli--Th\'era.
\end{df}
\begin{df}
Let $\phi$ be an OT function on $X$.
\begin{itemize}
\item[$(i)$] The \emph{OT ball} of $x\in X$ is:
\[
B_x^\phi:=\{y\in X\st d(x,y)\le-\phi(x,y)\}.
\]
If no confusion arises as to which OT function is considered, we will write $B_x$ in place of $B_x^\phi$. This gives rise to a ball space $(X,\mB^\phi)$, where
\[
\mB^\phi:=\{B_x\st x\in X\}.
\]
\item[$(ii)$] The \emph{OT ball space generated by an OT element }$x_0$ is $(B_{x_0},\mB^\phi_{x_0})$ where
\[
\mB^\phi_{x_0}:=\{B_x\st x\in B_{x_0}\}.
\]
\end{itemize}
\end{df}
\textbf{In this subsection, if an OT element $x_0\in X$ has been fixed, we will write for brevity $B_0$ in place of $B_{x_0}$.}\par
It is worth noting that if we are given a CK function $\varphi$, we may define $\phi$ by:
\begin{equation}
\label{smolphi}
\phi(x,y):=\varphi(y)-\varphi(x).
\end{equation}
The following fact is straightforward to prove.
\begin{fkt}
\label{smol}
If $\varphi$ is a CK function, then the function $\phi:X\times X\to\R$ defined in (\ref{smolphi}) is an OT function. Moreover, every $x\in X$ is an OT element for $\phi$.
\end{fkt}
As we know from Corollary \ref{co00}, the CK ball space is strongly contractive. A similar result can be shown for the OT ball space.
\begin{lem}
\label{lem0.2}
Take a metric space $\xD$ and $\phi:X\times X\to\R$ a function satisfying (b) and (c) in Definition \ref{whoa}. Then the following assertions hold, for every $x\in X$.
\begin{enumerate}
\item $x\in B_x$.
\item If $y\in B_x$ then $B_y\subseteq B_x$.
\item If $y\in B_x\setminus\{x\}$ then $B_y\subsetneq B_x$ and $\phi(x,y)<\phi(y,x)$.
\end{enumerate}
\end{lem}
\pruf 
$ $\\
(1): Indeed, $d(x,x)=-\phi(x,x)=0$.\\
(2): Take $y\in B_x$, i.e.,
\[
d(x,y)\le-\phi(x,y).\]
Take any $z\in B_y$, then
\[
d(y,z)\le-\phi(y,z).\]
By condition (c) for an OT function we get
\[
d(x,z)\le d(x,y)+d(y,z)\le-\phi(x,y)-\phi(y,z)\le -\phi(x,z),\]
so $z\in B_x$ and, as a result, $B_y\subseteq B_x$.\\
(3): Let $y\in B_x$ and $y\ne x$. We wish to show that $x\notin B_y$. Suppose that $x\in B_y$. Then $d(y,x)\le-\phi(y,x)$ and by conditions (b) and (c) for an OT function we get
\[
0<d(y,x)+d(x,y)\le-\phi(y,x)-\phi(x,y)\le-\phi(y,y)=0,\]
contradiction. Thus $x\notin B_y$ and so $B_y\subsetneq B_x$. Clearly, this also implies
\[
-\phi(y,x)<d(x,y)\le-\phi(x,y).
\]
\pruuf
Lemma \ref{lem0.2} instantly yields the following corollary.
\begin{co}
\label{co0}
For an OT function $\phi$ on $X$, the ball space $\mB^\phi$ is strongly contractive. Furthermore, for a fixed OT element $x_0$ for $\phi$ in $X$ the OT ball space $\mB^\phi_{x_0}$ is also strongly contractive and all of its balls are contained in $B_0$.
\end{co}
As stated in Fact \ref{smol}, for the OT function $\phi$ defined in (\ref{smolphi}) every $x\in X$ is an OT element. While this doesn't have to be true in general for any OT function $\phi$, this property turns out to be `hereditary' in the following sense.
\begin{lem}
Let $\phi$ be an OT function. If $x_0\in X$ is an OT element for $\phi$ in $X$ and $x\in B_0$ then also $x$ is an OT element for $\phi$ in $X$.
\end{lem}
\pruf
Let $r\in\R$ be such that 
\[
\inf_{y\in X}\phi(x_0,y)\ge r.
\]
Take any $x\in B_0$. Note that $0\le d(x_0,x)\le-\phi(x_0,x)$. For every $y\in X$ we have 
\[
r\le\phi(x_0,y)\le\phi(x_0,x)+\phi(x,y),\]
so
\[
\phi(x,y)\ge r-\phi(x_0,x).\]
In particular,
\[
\inf_{y\in X}\phi(x,y)\ge r-\phi(x_0,x)\ge r.
\]
\pruuf
As stated in Proposition \ref{pr1}, there is an equivalence between completeness of a metric space and spherical completeness of the respective CK ball spaces. A similar result can be shown for the OT ball spaces. For that we will need to state an auxiliary lemma first.
\begin{lem}
\label{hlap1}
Let $(X,d)$ be a metric space, $\phi$ an OT function on $X$ and $x_0$ an OT element for $\phi$ in $X$. Moreover, let $\mN\subseteq\mB^\phi_{x_0}$ be a nest of balls and write $\mN=\{B_x\st x\in A\}$ for some set $A\subseteq B_0$. Then for every $x,y\in A$ we have
\begin{equation}\label{hlap2}
d(x,y)\le|\phi(x_0,x)-\phi(x_0,y)|.
\end{equation}
Moreover, the following statements are equivalent for every $x,y\in A$:
\begin{itemize}
\item[$(i)$] $y\in B_x$,
\item[$(ii)$] $\phi(x,y)\le\phi(y,x)$,
\item[$(iii)$] $\phi(x_0,y)\le\phi(x_0,x)$.
\end{itemize}
\end{lem}
\pruf
For every $x,y\in A$ either $x\in B_y$ or $y\in B_x$ since $\mN$ is a nest, so
\begin{equation}\label{hlap3}
d(x,y)\le\max\{-\phi(x,y),-\phi(y,x)\}.
\end{equation}
If the above maximum is equal to $-\phi(x,y)$, we have
\[
d(x,y)\le-\phi(x,y)\le\phi(x_0,x)-\phi(x_0,y)\le|\phi(x_0,x)-\phi(x_0,y)|.\]
Similarly, if the maximum is equal to $-\phi(y,x)$, we have
\[
d(x,y)\le-\phi(y,x)\le\phi(x_0,y)-\phi(x_0,x)\le|\phi(x_0,x)-\phi(x_0,y)|.\]
Either way, we deduce (\ref{hlap2}).\par
To prove $(i)\Leftrightarrow(ii)$ assume $y\in B_x$. If $y=x$ then $(ii)$ is trivial. If $y\ne x$ then, by assertion (3) of Lemma \ref{lem0.2} we have 
\[
-\phi(y,x)<-\phi(x,y).
\]  
Hence $(ii)$ follows. Conversely, if $y\notin B_x$ (in particular, $y\ne x$) then $x\in B_y\setminus\{y\}$. As a result, again by assertion (3) of Lemma \ref{lem0.2}, $-\phi(x,y)<-\phi(y,x)$.\par 
To prove $(i)\Leftrightarrow(iii)$ we proceed as follows. If $y\in B_x$ then 
\[
0\le d(x,y)\le-\phi(x,y)\le-\phi(x_0,y)+\phi(x_0,x),
\] 
thus $\phi(x_0,x)\ge\phi(x_0,y)$. For the converse, if $y\notin B_x$ then $x\in B_y$ and so
\[
0<d(x,y)\le-\phi(y,x)\le-\phi(x_0,x)+\phi(x_0,y),
\] 
hence $\phi(x_0,x)<\phi(x_0,y)$.
\pruuf
\begin{pr}
\label{pr1.2}
A metric space $\xD$ is complete if and only if the OT ball space $(B^\phi_{x_0},\mB^\phi_{x_0})$ is spherically complete for every OT function $\phi$ on $X$ and every OT element $x_0$ for $\phi$ in $X$.
\end{pr}
\pruf
Suppose that for every OT function $\phi$ and every OT element $x_0$ for $\phi$ in $X$ the ball space $(B_0,\mB^\phi_{x_0})$ is spherically complete. We wish to show that the ball space $(X,\mB^\varphi)$ is spherically complete for every CK function $\varphi$ on $X$, which by Proposition \ref{pr1} will yield the completeness of the space $X$.\par
Take any CK function $\varphi$ on $X$, consider the ball space $(X,\mB^\varphi)$ and fix any nest of balls $\mN$ in $\mB^\varphi$. Pick some $B^\varphi_{x_0}\in\mN$ and consider the nest 
\[
\mN_0=\{ B\in\mN\st B\subseteq B^\varphi_{x_0}\}.
\]
By Fact \ref{smol} $x_0$ is an OT element for the OT function $\phi$ defined as in (\ref{smolphi}). Moreover, for every $x\in X$ such that $B^\varphi_x\subseteq B^\varphi_{x_0}$, we have $B^\varphi_{x}=B^\phi_{x}$, hence $\mN_0$ is a nest in the OT ball space $(B_{x_0}^\phi,\mB_{x_0}^\phi)$. By assumption we then obtain that $\emptyset\neq\bigcap\mN_0=\bigcap\mN$.\par
For the converse, assume that $X$ is complete. Fix any OT function $\phi$ on $X$ and any OT element $x_0$ for $\phi$ in $X$. Take a nest of balls $\mN$ in the ball space $\mB^\phi_{x_0}$ and write $\mN=\{B_x\st x\in A\}$ for some set $A\subseteq B_0$. By assumption on $x_0$ there exists
\[
r:=\inf_{x\in A}\phi(x_0,x)\in\R.\]
Let $(x_n)_{n\in\N}$ be a sequence of elements in $A$ such that
\[
\lim_{n\to\infty}\phi(x_0,x_n)=r.\]
Then $(\phi(x_0,x_n))_{n\in\N}$ is a Cauchy sequence because it converges to $r$. By Lemma \ref{hlap1} the sequence $(x_n)_{n\in\N}$ is also Cauchy. Since $X$ is complete, $(x_n)_{n\in\N}$ converges to an element $z\in X$. We want to show that $z\in\bigcap\mN$ or, equivalently, that $z\in B_x$ for every $x\in A$. Fix an arbitrary element $x\in A$. If $\phi(x_0,x)=r$ (in particular, the infimum is achieved) then by Lemma \ref{hlap1} we get that $x=z$, because
\[
d(x_n,x)\le|\phi(x_0,x_n)-\phi(x_0,x)|=|\phi(x_0,x_n)-r|\to0,\]
showing that $x$ is a limit of $(x_n)_{n\in\N}$. Hence in this case we obtain that $z\in B_x$ trivially. Therefore we may assume that $\phi(x_0,x)>r$. Then from the definition of $(x_n)_{n\in\N}$ we obtain the existence of $N\in\N$ such that, for every $n\ge N$, we have $\phi(x_0,x_n)\le\phi(x_0,x)$. This, by Lemma \ref{hlap1}, is equivalent to $\phi(x,x_n)\le\phi(x_n,x)$. Therefore, for every $n\ge N$,
\[
d(x,x_n)\le\max\{-\phi(x,x_n),-\phi(x_n,x)\}=-\phi(x,x_n),
\]
where the first inequality is deduced similarly to (\ref{hlap3}). Taking $\limsup$ on both sides we get
\[
d(x,z)\le \limsup_{n\to\infty}-\phi(x,x_n)\le -\phi(x,z),
\] 
so that $z\in B_x$. Since $x\in A$ was an arbitrary element, we get $z\in\bigcap\mN$ as claimed.
\pruuf
\begin{re}
Proposition \ref{pr1.2} does in general not hold for the ball space $\mB^\phi$ in place of $\mB^\phi_{x_0}$. Take the complete metric space $\R$, where $d(x,y)=|x-y|$, and consider $\phi:\R\times\R\to\R$ defined as:
\[
\phi(x,y)=\begin{cases} x-y&\text{if }x\ne0\\ 0&\text{if }x=0\end{cases}.
\]
This is an OT function and yields balls of the form
\[
B_x=\{y\in X\st |x-y|\le y-x\}=[x,\infty).
\]
In the corresponding ball space we have a nest of balls with empty intersection, namely,
\[
\{[n,+\infty)\st n\in\N\}.
\]
\end{re}
Armed with the theory introduced so far, we can prove an important property of OT (and as a result, also CK) ball spaces in a complete metric space.
\begin{pr}
\label{tw1.1}
Let $\xD$ be a complete metric space.
\begin{enumerate}
\item If $\phi$ is an OT function on $X$ then for every OT element $x_0$ for $\phi$ in $X$ there exists an element $a\in B_0$ such that $B^\phi_a=\{a\}$.
\item If $\varphi$ is a CK function on $X$ then for every $x\in X$ there exists $a\in B^\varphi_x$ such that $B^\varphi_a=\{a\}$.
\end{enumerate}
\end{pr}
\pruf
Assertion (2) follows from assertion (1) by Fact \ref{smol}. To prove assertion (1) let $\phi$, $x_0$ and $B_0$ be as in the assumption of the Proposition. By Proposition \ref{pr1.2} the OT ball space $\mB^\phi_{x_0}$ is spherically complete, and by Corollary \ref{co0} it is strongly contractive. Theorem \ref{poetycko} yields the result.
\pruuf

\subsection{Generalized Caristi--Kirk ball spaces}
\label{CKBSG}
Consider a function $\tvarphi:X\to(-\infty,+\infty]$ which is lower semicontinuous, bounded from below and not constantly equal to $+\infty$. We will call such $\tvarphi$ a \emph{CK}$^{\infty}$ \emph{function on} $X$. In this setting we may define the \emph{CK}$^\infty$ \emph{balls} as follows:
\[
B_x^{\tvarphi}:=\{x\in X\st \tvarphi(y)+d(x,y)\le\tvarphi(x)\}.
\]
If an element $x_0\in X$ satisfies $\tvarphi(x)<+\infty$, we will call it a \emph{CK element for }$\tvarphi$\emph{ in }$X$ (or simply a \emph{CK element}).\par
An easy observation is that every CK function is a CK$^\infty$ function. However, for a CK$^\infty$ function $\tvarphi$, setting
\begin{equation}
\label{smoltphi}
\phi(x,y):=\tvarphi(y)-\tvarphi(x)
\end{equation}
as we did in the CK case (\ref{smolphi}), may not make sense. Indeed, in the case $\tvarphi(x)=+\infty=\tvarphi(y)$ there is no natural choice for the value of $\phi(x,y)$.\par
\textbf{In this subsection, if a CK element $x_0$ is fixed, we will write $B_0$ in place of $B^\tvarphi_{x_0}$.}\par
For a CK element $x_0$ we define the CK$^\infty$ \emph{ball space generated by }$x_0$ as the ball space $(B_0,\mB^{\tvarphi}_{x_0})$, where:
\[
\mB^{\tvarphi}_{x_0}:=\{B^\tvarphi_x\st x\in B_0\}.
\]
Note that in general the ball space $\{B_x^\tvarphi\st x\in X\}$ is not strongly contractive. Indeed, if $x,y\in X$, $x\ne y$, satisfy $\tvarphi(x)=\tvarphi(y)=+\infty$, then $B^\tvarphi_x=X=B^\tvarphi_y$. However, if we work inside a CK$^\infty$ ball space, strong contractiveness holds, as stated in the following lemma. 
\begin{lem}
\label{lem0.3}
Take a metric space $\xD$, any function $\tvarphi:X\to(-\infty,+\infty]$. Then the following assertions hold for every $x\in X$.
\begin{enumerate}
\item $x\in B^\tvarphi_x$.
\item If $y\in B^\tvarphi_x$ then $B^\tvarphi_y\subseteq B^\tvarphi_x$ and $\tvarphi(y)\le\tvarphi(x)$.
\item Let $x\in X$ be such that $\tvarphi(x)<+\infty$ and let $y\in B^\tvarphi_x\setminus\{x\}$. Then $B^\tvarphi_y\subsetneq B^\tvarphi_x$ and $\tvarphi(y)<+\infty$.
\end{enumerate}
\end{lem}
\pruf$ $\\
(1): Indeed, $\tvarphi(x)+d(x,x)=\tvarphi(x)$.\\
(2): If $\varphi(x)=+\infty$ then $B^\tvarphi_x=X$ and $B^\tvarphi_y\subseteq B^\tvarphi_x$ as well as $\tvarphi(y)\le\tvarphi(x)$ trivially. Now assume that $\tvarphi(x)<+\infty$ and $y\in B^\tvarphi_x$. Then also $\tvarphi(y)<+\infty$ because
\[
\tvarphi(y)\le \tvarphi(x)-d(x,y)\le\tvarphi(x)<+\infty.
\]
Take any $z\in B^\tvarphi_y$. Through the same reasoning as above, we can see that $\tvarphi(z)<+\infty$ and therefore we may write
\begin{align*}
d(x,z)&\le d(x,y)+d(y,z)\\
&\le\tvarphi(x)-\tvarphi(y)+\tvarphi(y)-\tvarphi(z)\\
&=\tvarphi(x)-\tvarphi(z),
\end{align*}
which shows that $z\in B^\tvarphi_x$. Since $z\in B^\tvarphi_y$ was arbitrary, we deduce that $B^\tvarphi_y\subseteq B^\tvarphi_x$.\\
(3): Assume that $\tvarphi(x)<+\infty$ and $y\in B^\tvarphi_x\setminus\{x\}$. From assertion (2) of our lemma we know that $\tvarphi(y)<+\infty$ and $B^\tvarphi_y\subseteq B^\tvarphi_x$. Suppose that $x\in B^\tvarphi_y$. In that case
\[
\tvarphi(y)+d(x,y)\le\tvarphi(x)
\]
and
\[
\tvarphi(x)+d(x,y)\le\tvarphi(y).
\]
Adding up these inequalities, taking into account that $\tvarphi(x)<+\infty$ and $\tvarphi(y)<+\infty$, we obtain
\[
0<2d(x,y)\le \tvarphi(x)-\tvarphi(y)+\tvarphi(y)-\tvarphi(x)=0,
\]
contradiction. So we must have $x\notin B^\tvarphi_y$, hence $B^\tvarphi_y\subsetneq B^\tvarphi_x$.\qedhere
\pruuf
From assertion (3) of Lemma \ref{lem0.3}, we obtain:
\begin{co}
\label{hlp}
Let $x_0$ be a CK element for a CK$^\infty$ function $\tvarphi$. Then for every $y\in B_0$, $y$ is also a CK element for $\tvarphi$. Further, $\tvarphi |_{B_0}$ is a CK function and $(B_0,\mB^\tvarphi_{x_0})$ is a CK ball space in the sense of Section \ref{CKBS}.
\end{co}

Before we state another property of CK$^\infty$ balls, it is worth noting that the proofs of Lemma \ref{lem0.3} and Lemma \ref{closedbols} are similar (or, at times, identical) to the original proof of Lemma \ref{lem0}, which can be found in \cite{karticle}. This comes from the fact that for a CK element $x_0$, the set $\tvarphi(B_0)$ does not contain infinity, so these balls keep the properties of the `original' CK balls.
\begin{lem}
\label{closedbols}
For every $x\in X$ and every CK$^\infty$ function $\tvarphi$, the ball $B^\tvarphi_x$ is closed in the topology induced by the metric.
\end{lem}
\pruf
The complement $\{y\in X\st d(x,y)+\tvarphi(y)>\tvarphi(x)\}$ of $B^\tvarphi_x$ is the preimage of the final segment $(\tvarphi(x),+\infty]$ of $(-\infty,+\infty]$, which is open in the Scott topology, under the function
\[
X\ni y\stackrel{\psi}{\longmapsto} d(x,y)+\tvarphi(y).
\]
Whenever $\tvarphi$ is lower semicontinuous, then so is $\psi$ and this preimage is open in $X$.
\pruuf
We are now ready to prove a result analogous to Propositions \ref{pr1} and \ref{tw1.1}.
\begin{pr}
\label{tw1.3}
Let $(X,d)$ be a complete metric space and $\tvarphi$ be a CK$^\infty$ function on $X$. If $x_0\in X$ then there exists $a\in B^\tvarphi_{x_0}$ such that $B^\tvarphi_a=\{a\}$.
\end{pr}
\pruf
Consider a complete metric space $(X,d)$, fix any element $x_0\in X$ and consider the ball $B_0:=B_{x_0}^\tvarphi$.\par
Assume first that $x_0$ is a CK element for $\tvarphi$ in $X$. As we know from Lemma \ref{closedbols}, $B_0$ is closed, hence complete. Moreover, the function $\varphi:=\tvarphi|_{B_0}$ is a CK function.
Note that for every $x\in B_0$ we have
\begin{align*}
B_x^\varphi&=\{y\in B_0\st d(x,y)\le\varphi(x)-\varphi(y)\}\\
&\subseteq B_x^\tvarphi=\{y\in X\st\tvarphi(y)+d(x,y)\le\tvarphi(x)\}.
\end{align*}
We wish to show that the above sets are equal. By assertion (2) of Lemma~\ref{lem0.3} we know that $B_x^\tvarphi\subseteq B_{0}$. On $B_0$ we have $\varphi=\tvarphi$ so that the values of $\tvarphi$ are finite and $\tvarphi(y)+d(x,y)\le\tvarphi(x)$ is equivalent to $d(x,y)\le\varphi(x)-\varphi(y)$. This yields $B_x^\tvarphi\subseteq B_x^\varphi$.\par
Since $\varphi$ is a CK function on a complete metric space $B_{0}$, we may apply assertion (2) of Proposition \ref{tw1.1} to the CK ball space $(B_{0},\mB^\varphi_{x_0})$, to acquire $a\in B_{0}$ such that
\[
\{a\}=B_a^\varphi=B_a^\tvarphi.
\]
Assume now that $x_0\in X$ is not a CK element for $\tvarphi$. Then we obtain that $B_{0}=X$. Inside the ball $B_{0}$ we may thus find a CK element $x_1$ for $\tvarphi$. From what we have proved above there exists $a\in B^\tvarphi_{x_1}\subseteq X=B_{0}$ such that $B^\tvarphi_a=\{a\}$.
\pruuf

\section{Applications of Proposition \ref{tw1.1}}\label{applications}
In this section we give simple proofs for a number of known theorems, in versions that involve OT functions, by applying Proposition \ref{tw1.1}. Note that the multivalued Caristi--Kirk FPT, Ekeland's Principle and Takahashi's Theorem have already been proved in the OT form in \cite{ot} using the Oettli--Th\'era Theorem. The original versions of these theorems are listed in Section~\ref{original}.
\begin{tw}[Caristi--Kirk FPT, OT form]
Let $\xD$ be a complete metric space and $\phi$ an OT function on $X$. If a function $f:X\to X$ satisfies:
\begin{equation}
\label{ck1.2}
\forol{x\in X}d(x,f(x))\le-\phi(x,f(x)),
\end{equation}
then $f$ has a fixed point on $X$, i.e., there exists an element $a\in X$ such that $f(a)=a$.
\end{tw}
\pruf
Condition (\ref{ck1.2}) implies that for every $x\in X$ we have
\[
f(x)\in B_x.
\]
Proposition \ref{tw1.1} gives us the existence of $a\in X$ such that $B_a=\{a\}$. In particular, since $f(a)\in B_a$, we have $f(a)=a$.
\pruuf
\begin{tw}[Caristi--Kirk FPT, multivalued version, OT form]
Let $\xD$ be a complete metric space and $\phi$ an OT function on $X$.\par If a function $F:X\to\mP(X)$ satisfies:
\begin{equation}
\label{ck2.2}
\forol{x\in X}\egzisc{y\in F(x)}d(x,y)\le-\phi(x,y),
\end{equation}
then $F$ has a fixed point on $X$, i.e., there exists $a\in X$ such that $a\in F(a)$.
\end{tw}
\pruf
Condition (\ref{ck2.2}) means that for every $x\in X$ there exists $y\in F(x)\cap B_x$. In particular, for $x=a$ with $a\in X$ given by Proposition \ref{tw1.1} we obtain
\[
y\in F(a)\cap B_a\subseteq\{a\},
\]
whence $a=y\in F(a)$.
\pruuf
\begin{tw}[Basic Ekeland's Principle, OT form]
Let $\xD$ be a complete metric space and $\phi$ an OT function on $X$. There exists $a\in X$ such that
	\begin{equation}
	\label{CC3}
	\forol{x\in X\setminus\{a\}}-\phi(a,x)<d(a,x).
	\end{equation}
\end{tw}
\pruf
Property (\ref{CC3}) is equivalent to $B_a=\{a\}$ and the existence of such $a\in X$ follows from Proposition \ref{tw1.1}.
\pruuf
\begin{tw}[Altered Ekeland's Principle, OT form]
\label{aep}
Let $\xD$ be a complete metric space and $\phi$ an OT function on $X$. For any $\gamma>0$ and any OT element $x_0$ for $\phi$ in $X$ there exists $a\in X$ such that
\begin{equation}\label{CC4}
\forol{x\in X\setminus\{a\}} -\phi(a,x)<\gamma d(a,x)
\end{equation}
and
\begin{equation}\label{CC5}
-\phi(x_0,a)\ge\gamma d(x_0,a).
\end{equation}
\end{tw}
\pruf
Since $\gamma>0$, the function $\psi:=\gamma^{-1}\phi$ is an OT function on $X$, so we can work with $\psi$ and the respective ball space $\mB^\psi$.\par
We apply Proposition \ref{tw1.1} to the given complete metric space $X$, the function $\psi$ and the ball $B_0:=B_{x_0}^\psi$. This gives us the existence of an element $a\in B_0$ such that $B_a^\psi=\{a\}$. Now, the assertion $a\in B_0$ means that 
\[
d(x_0,a)\le-\psi(x_0,a)=-\gamma^{-1}\phi(x_0,a),
\]
which is equivalent to property (\ref{CC5}). Similarly, $B_a^\psi=\{a\}$ implies
\[
\forol{x\in X\setminus\{a\}} d(a,x)>-\psi(a,x)=-\gamma^{-1}\phi(a,x),
\]
which is equivalent to property (\ref{CC4}). 
\pruuf
\begin{tw}[Ekeland's Usual Variational Theorem, OT form]
Let $\xD$ be a complete metric space and $\phi$ an OT function on $X$. Fix $\varepsilon\ge0$ and $x_0\in X$ such that $-\varepsilon\leq \inf_{x\in X}\phi(x_0,x)$. Then for any $\gamma>0$ and $\delta\ge0$ with $\gamma\delta\ge\varepsilon$ there exists $a\in X$ such that $d(a,x_0)\le\delta$ and $a$ is the strict minimum point of the function $\phi_\gamma:X\to(-\infty,+\infty]$ defined as
\[
\phi_\gamma(x)=\phi(a,x)+\gamma d(x,a).
\]
\end{tw}
\pruf
Take $\varepsilon\ge0$ and $x_0$ as in the assumptions of the theorem, and fix arbitrary real numbers $\gamma>0$ and $\delta\ge0$ such that $\gamma\delta\ge\varepsilon$. The function $\psi:=\gamma^{-1}\phi:X\times X\to(-\infty,+\infty]$ is an OT function on $X$, so we can apply Proposition \ref{tw1.1} with the function $\psi$ and $B_0:=B^\psi_{x_0}$ (note that $x_0$ is an OT element for $\psi$ in $X$). We deduce the existence of $a\in B_0$ such that $B^\psi_a=\{a\}$. Now, the property $a\in B_0$ means that
\[
d(x_0,a)\le-\psi(x_0,a),
\]
or in other words:
\[
\gamma d(x_0,a)\le-\phi(x_0,a)\le-\inf_{x\in X}\phi(x_0,x)\le\varepsilon\le\gamma\delta.
\]
Thus,
\[
d(a,x_0)=d(x_0,a)\le\delta.
\]
The property $B_a=\{a\}$ means that for every $x\in X\setminus\{a\}$ we have that
\[
d(x,a)>-\psi(a,x)=-\gamma^{-1}\phi(a,x).
\]
From this we obtain that
\[
\phi_\gamma(x)=\phi(a,x)+\gamma d(x,a)>0=\phi_\gamma(a),
\]
which means that $a$ is the strict minimum point of the function $\phi_\gamma$.
\pruuf
\begin{df}
Let $\xD$ be a metric space. Take $\gamma\in(0,\infty)$ and $a,b\in X$. The \emph{petal associated with }$\gamma$\emph{ and }$a,b$ is the subset $P_\gamma(a,b)$ of $X$ defined as follows:
\[
P_\gamma(a,b)=\{y\in X\st\gamma d(y,a)+d(y,b)\le d(a,b)\}.
\]
\end{df}
\begin{tw}[Flower Petal Theorem]
Let $M$ be a complete subset of a metric space $\xD$. Take $x_0\in M$ and $b\in X\setminus M$. Then for each $\gamma>0$ there exists $a\in P_\gamma(x_0,b)\cap M$ such that
\[
P_\gamma(a,b)\cap M=\{a\}.
\]
\end{tw}
\pruf
We use the notation from the assertion of the theorem. As $\gamma>0$, the function $\varphi:M\to\R$ given by
\[
\varphi(x):=\gamma^{-1}d(x,b)
\]
is a CK function on $M$. In this setting we have, for every $x\in M$,
\[
P_\gamma(x,b)\cap M=\{y\in M\st d(x,y)\le\varphi(x)-\varphi(y)\}=B^\varphi_x.
\]
To conclude we use assertion (2) of Proposition \ref{tw1.1} with $M$ in place of $X$ and $x:=x_0$, which yields the existence of $a\in B_{x_0}^\varphi=P_\gamma(x_0,b)\cap M$ such that
\[
\{a\}=B^\varphi_a=P_\gamma(a,b)\cap M.
\]
\pruuf
\begin{tw}[Takahashi, OT form]
Let $\xD$ be a complete metric space, $\phi$ an OT function on $X$ and $x_0\in X$ an OT element for $\phi$ in $X$. Assume that for every $u\in B_{x_0}$ with $\inf_{x\in X}\phi(u,x)<0$ there exists $v\in X$ such that $v\ne u$ and $d(u,v)\le-\phi(u,v)$. Then there exists $a\in B_{x_0}$ such that $\inf_{x\in X}\phi(a,x)=0$.
\end{tw}
\pruf
Proposition \ref{tw1.1} gives us the existence of $a\in B_{x_0}$ such that $B_a=\{a\}$. If $\inf_{x\in X}\phi(a,x)<0$, then by assumption there would exist $v\in X\setminus\{a\}$ such that $d(a,v)\le-\phi(a,v)$, which would mean that $B_a$ is not a singleton, contradiction. So $\inf_{x\in X}\phi(a,x)\geq 0$, but $\phi(a,a)=0$ which proves the claim.
\pruuf
\begin{tw}[Oettli-Th\'era]
Let $\xD$ be a complete metric space, $\phi$ an OT function on $X$ and $x_0\in X$ an OT element for $\phi$ in $X$. Let $\Psi\subseteq X$ have the property that
\begin{equation}
\forol{x\in B_{x_0}\setminus\Psi}\egzisc{y\in X\setminus\{x\}}d(x,y)\le-\phi(x,y).
\end{equation}
Then there exists $a\in B_{x_0}\cap\Psi$.
\end{tw}
\pruf
From Proposition \ref{tw1.1} there exists $a\in B_{x_0}$ such that $B_a=\{a\}$. If $a\notin \Psi$ then, by assumption, $B_a$ would contain another element $y\ne a$, which would mean that $B_a$ is not a singleton, contradiction.
\pruuf
\section{Applications of Proposition \ref{tw1.3}}\label{original}
Many of the theorems mentioned in the previous section have been originally stated and proved using the CK function $\varphi$. By Fact \ref{smol}, proving the version involving $\phi$, through (\ref{smolphi}), will also automatically prove the version involving $\varphi$. However, many sources (e.g., \cite{penot}, \cite{araya}) cite the theorems in a CK$^\infty$ form. As already remarked, we cannot directly define an OT function from a CK$^\infty$ function. Nevertheless, we can use Proposition \ref{tw1.3} to prove these versions in the same way we did in the previous section using Proposition \ref{tw1.1}. Since the proofs are analogous to the ones stated in Section \ref{applications}, we will leave them to the reader.\par
Note that here we do not include the Oettli-Th\'era Theorem (since it has originally been stated in the OT form) nor the Flower Petal Theorem (since it does not include either of the functions).\par
For the following theorems, fix a complete metric space $\xD$ and a CK$^\infty$ function $\tvarphi$ on $X$.
\begin{tw}[Caristi-Kirk FPT, CK$^\infty$ form]
If a function $f:X\to X$ satisfies
\[
\forol{x\in X}\tvarphi(f(x))+d(x,f(x))\le\tvarphi(x)
\]
then $f$ has a fixed point on $X$, i.e., there exists an element $a\in X$ such that $f(a)=a$.
\end{tw}
\begin{tw}[Caristi-Kirk FPT, multivalued version, CK$^\infty$ form]
If a function $F:X\to\mP(X)$ satisfies
\begin{equation}
\label{ckmv}
\forol{x\in X}\egzisc{y\in F(x)}\tvarphi(y)+d(x,y)\le\tvarphi(x)
\end{equation}
then $F$ has a fixed point on $X$, i.e., there exists $a\in X$ such that $a\in F(a)$.
\end{tw}
\begin{tw}[Basic Ekeland's Principle, CK$^\infty$ form]
There exists $a\in X$ such that
\[
\forol{x\in X\setminus\{a\}}\tvarphi(a)<\tvarphi(x)+d(a,x).
\]
\end{tw}
\begin{tw}[Altered Ekeland's Principle, CK$^\infty$ form]
For all $\gamma>0$ and any $x_0\in X$ there exists $a\in X$ such that
\[
\forol{x\in X\setminus\{ a\}}\tvarphi(a)<\tvarphi(x)+\gamma d(a,x)
\]
and
\[
\tvarphi(a)\le\tvarphi(x_0)-\gamma d(a,x_0).
\]
\end{tw}
\begin{tw}[Ekeland's Usual Variational Theorem, CK$^\infty$ form]
Let $\varepsilon\ge0$ and $x_0\in X$ be such that $\tvarphi(x_0)\le\inf\tvarphi(X)+\varepsilon$. Then for any $\gamma>0$ and $\delta\ge0$ with $\gamma\delta\ge\varepsilon$ there exists $a\in X$ such that $d(a,x_0)\le\delta$ and $a$ is the strict minimum point of the function
\[
\tvarphi_\gamma(x)=\tvarphi(x)+\gamma d(x,a).
\]
\end{tw}
\begin{tw}[Takahashi, CK$^\infty$ form]
Suppose that for each $u\in X$ with $\inf_{x\in X}\tvarphi(x)<\tvarphi(u)$ there exists $v\in X$ such that $v\ne u$ and $\tvarphi(v)+d(u,v)\le\tvarphi(u)$. Then there exists $a\in X$ such that $\inf_{x\in X}\tvarphi(x)=\tvarphi(a)$.
\end{tw}

\newpage
\bibliographystyle{siam}
\bibliography{refs2}

\begin{thebibliography}{10}

\bibitem{araya}
{\sc Y.~Araya}, {\em On generalizing {T}akahashi's nonconvex minimization
  theorem}, Appl. Math. Lett., 22 (2009), pp.~501--504.

\bibitem{caristi}
{\sc J.~Caristi}, {\em Fixed point theorems for mappings satisfying inwardness
  conditions}, Trans. Amer. Math. Soc., 215 (1976), pp.~241--251.

\bibitem{ekeland}
{\sc I.~Ekeland}, {\em On the variational principle}, J. Math. Anal. Appl., 47
  (1974), pp.~324--353.

\bibitem{kirk}
{\sc W.~A. Kirk}, {\em Caristi's fixed point theorem and metric convexity},
  Colloq. Math., 36 (1976), pp.~81--86.

\bibitem{kbook}
{\sc F.-V. Kuhlmann and K.~Kuhlmann}, {\em Ball spaces -- a basic framework for
  fixed point theorems: ball spaces and spherical completeness, in
  preparation}.

\bibitem{kbonus}
\leavevmode\vrule height 2pt depth -1.6pt width 23pt, {\em A common
  generalization of metric and ultrametric fixed point theorems}, Forum Math.,
  27 (2015), pp.~303--327.

\bibitem{kbonus2}
\leavevmode\vrule height 2pt depth -1.6pt width 23pt, {\em Correction to ``a
  common generalization of metric, ultrametric and topological fixed point
  theorems''}, Forum Math., 27 (2015), pp.~329--330.

\bibitem{karticle}
{\sc F.-V. Kuhlmann, K.~Kuhlmann, and M.~Paulsen}, {\em The {C}aristi-{K}irk
  fixed point theorem from the point of view of ball spaces}, J. Fixed Point
  Theory Appl., 20 (2018), pp.~Art. 107, 9.

\bibitem{meghea}
{\sc I.~Meghea}, {\em Ekeland variational principle}, Old City Publishing,
  Philadelphia, PA; \'{E}ditions des Archives Contemporaires, Paris, 2009.

\bibitem{ot}
{\sc W.~Oettli and M.~Th\'{e}ra}, {\em Equivalents of {E}keland's principle},
  Bull. Austral. Math. Soc., 48 (1993), pp.~385--392.

\bibitem{penot}
{\sc J.-P. Penot}, {\em The drop theorem, the petal theorem and {E}keland's
  variational principle}, Nonlinear Anal., 10 (1986), pp.~813--822.

\bibitem{takahashi}
{\sc W.~Takahashi}, {\em Existence theorems generalizing fixed point theorems
  for multivalued mappings}, in Fixed point theory and applications
  ({M}arseille, 1989), vol.~252 of Pitman Res. Notes Math. Ser., Longman Sci.
  Tech., Harlow, 1991, pp.~397--406.

\end{thebibliography}
\end{document}